\documentclass[10pt]{amsart}
\usepackage{epic,eepic,latexsym}

 \newlength{\baseunit}               
 \newcount{\numlines}                
\newcommand{\point}{\vspace{3mm}\par\refstepcounter{subsection}{\bf \thesubsection.} }
\newcommand{\tpoint}[1]{\vspace{3mm}\par\refstepcounter{subsection}{\bf \thesubsection.} 
  {\em #1. ---} }
\newcommand{\epoint}[1]{\vspace{3mm}\par\refstepcounter{subsection}{\bf \thesubsection.} 
  {\em #1.} }
\newcommand{\bpoint}[1]{\vspace{3mm}\par\refstepcounter{subsection}{\bf \thesubsection.} 
  {\bf #1.} }

\newcommand{\bpf}{\noindent {\em Proof.  }}
\newcommand{\epf}{\qed \vspace{+10pt}}

\newcommand{\zed}{\mathbb{Z}}

\newcommand{\Q}{\mathbb{Q}}
\newcommand{\com}{\mathbb{C}}
\newcommand{\eff}{\mathbb{F}}

\newcommand{\proj}{\mathbb P}
\newcommand{\oh}{{\mathcal{O}}}
\newcommand{\cL}{{\mathcal{L}}}
\newcommand{\cm}{{\mathcal{M}}}
\newcommand{\ch}{{\mathcal{H}}}
\newcommand{\tC}{\tilde{C}}

\newcommand{\De}{\Delta}
\newcommand{\si}{\sigma}

\newcommand{\odd}{\operatorname{odd}}

\newcommand{\even}{\operatorname{even}}

\newcommand{\Pic}{\operatorname{Pic}}

\newcommand{\Bl}{\operatorname{Bl}}
\newcommand{\Sym}{\operatorname{Sym}}
\newcommand{\Aut}{\operatorname{Aut}}
\newcommand{\AAA}{\mathcal{A}}
\newcommand{\BBB}{\mathcal{B}}
\newcommand{\CCC}{\mathcal{C}}
\newcommand{\DDD}{\mathcal{D}}
\newcommand{\EEE}{\mathcal{E}}
\newcommand{\FFF}{\mathcal{F}}
\newcommand{\ZZZ}{\mathcal{Z}}
\newcommand{\grass}{\mathbb{G}}
\newcommand{\cited}{}

\begin{document}
\pagestyle{plain}
\title{Twelve points on the projective line, branched covers, 
and rational elliptic fibrations\footnote{1991 Mathematics Subject Classification:  Primary 14H10, Secondary 14H45}}
\author{Ravi Vakil}
\address{Dept. of Mathematics, MIT, Cambridge MA~02139}
\email{vakil@math.mit.edu}
\thanks{Partially supported by NSF Grant DMS--9970101}
\date{October 1, 1999.}
\begin{abstract}
The following divisors in the space $\Sym^{12} \proj^1$ of twelve
points on $\proj^1$ are actually the same: $(\AAA)$ the
possible locus of the twelve nodal fibers in a rational elliptic
fibration (i.e. a pencil of plane cubic curves); $(\BBB)$ degree 12
binary forms that can be expressed as a cube plus a square; $(\CCC)$
the locus of the twelve tangents to a smooth plane quartic from a
general point of the plane; $(\DDD)$ the branch locus of a degree 4
map from a hyperelliptic genus 3 curve to $\proj^1$;
$(\EEE)$ the branch locus of a degree 3 map from a genus 4 curve to
$\proj^1$ induced by a theta-characteristic; and several more.

The corresponding moduli spaces are smooth, but they are not all
isomorphic; some are finite \'{e}tale covers of others.  We describe
the web of interconnections between these spaces, and give monodromy,
rationality, and Prym-related consequences.  Enumerative consequences
include: (i) the degree of this locus is 3762 (e.g. there are 3762
rational elliptic fibrations with nodes above 11 given general points of
the base); (ii) if $C \rightarrow \proj^1$ is a cover as in $(\DDD)$,
then there are 135 different such covers branched at the same points;
(iii) the general set of 12 tangent lines that arise in $(\CCC)$ turn up
in 120 essentially different ways.

Some parts of this story are well-known, and some other parts were
known classically (to Zeuthen, Zariski, Coble, Mumford, and others).
The unified picture is surprisingly intricate and connects many
beautiful constructions, including Recillas' trigonal construction and
Shioda's $E_8$-Mordell-Weil lattice. 
\end{abstract}
\maketitle
\tableofcontents

{\parskip=12pt 

\section{Introduction}

In this paper, we explore the links between various moduli spaces
defined in Section \ref{dr} (and sketched in the abstract).  Although
some of the connections are well-known or classical, for completeness
of exposition we describe them in detail; some classical references
are discussed in Section \ref{cr}.

\epoint{Outline of the paper}
In Section \ref{pluto}, we interpret points of $\CCC$, $\DDD$, and
$\EEE$ as 2-torsion information on rational elliptic fibrations
(corresponding to points of $\AAA$); the Mordell-Weil lattice plays a
central role.  In Section \ref{uranus}, we relate $\CCC$, $\DDD$, and
$\EEE$ using Recillas' trigonal construction, and interpret points of
$\CCC$ and $\DDD$ as theta-characteristics of a genus 4 curve
corresponding to a point of $\EEE$.  As consequences, we give
monodromy and Prym-related results, including a theorem of
Mumford on hyperelliptic Pryms.  In the short Section \ref{ganymede},
we note that some of the connections are already visible in the
discriminant of a quartic or cubic polynomial.  In Section
\ref{saturn}, we study the actual locus of 12 points on the projective
line, the divisor $\ZZZ \subset \Sym^{12} \proj^1 \cong \proj^{12}$,
and compute its degree to be 3762.  In Section \ref{fq}, we suggest
further questions about two-torsion of elliptic fibrations.

These results can also be seen as related to Persson's enumeration of the
possible singular fibers of rational elliptic surfaces (\cite{p},
\cite{mi}); he studies how the singular fibers can degenerate, while
we study the locus of possible configurations.  (Note that \cite{p}
makes use of some of the standard constructions below.)

Many of the constructions in this
paper are reminiscent of \cite{b} Ch. VI, but there does not appear to
be a precise connection.

\epoint{Acknowledgements} 
The advice of I. Dolgachev has been invaluable.  In particular, 
the authour is grateful to him for suggesting that $\CCC$,
$\DDD$ and $\EEE$ might be related via Recillas' construction, and
for pointing out the beautiful classical reference \cite{co}.  
The author would also like to thank D. Allcock, A. Bertram, 
A. J. de Jong, N. Elkies, J. Harris, B. Hassett, and R. Smith for helpful
conversations.

\section{Definitions and Results}
\label{dr}

For convenience, we work over the complex numbers.  (All results are
true over an algebraically closed field of characteristic greater than
5, except possibly Theorem \ref{degZ} in characteristics 11 and 19.)
All moduli spaces are Deligne-Mumford stacks unless
otherwise indicated.

We now describe various loci of twelve points on the projective line,
and the relationships between them.  The constructions are invariant
under $\Aut(\proj^1) = PGL(2)$, and there are natural quotients in
each case.  Often the quotient is the more ``natural'' object, but for
consistency, we will deal with the ``framed'' construction.  (If $X$
is abstractly isomorphic to $\proj^1$, a {\em frame} for $X$ is a
choice of isomorphism to $\proj^1$.)

Most of these constructions have natural compactifications (or at
least partial compactifications).  It would be interesting to
understand how these relationships specialize to the boundary; the
reader will note that they often specialize well to the divisor where
two of the twelve points come together.  Many of the results of
\cite{e} (about six points on $\proj^1$) can be understood in this
context by letting the points come together in six pairs (for example
the results relating to genus 4 curves with vanishing theta
characteristic, and facts relating to the $E_8$-lattice, including the
numbers 136 and 120).  Many more of the results of \cite{e} have a
similar flavor, but the connections are not clear.

\bpoint{Degree 1 Del Pezzo surfaces and rational
elliptic fibrations}
\label{ellfib}
We recall relevant facts relating degree 1 Del Pezzo surfaces,
rational elliptic fibrations, and pencils of plane cubics.  (The
fundamental sources are \cite{shioda}, especially Section 10,
and \cite{manin} Ch. IV.)  By genus 1 fibration, we will mean a flat family
of reduced genus 1 curves over a smooth curve, with smooth total space
and at worst multiplicative reduction.  By elliptic fibration, we will
mean a genus 1 fibration with a choice of section.

A degree 1 Del Pezzo surface $X$ has a pencil of sections of the
canonical sheaf, with one base point $p_0$.  Then the pencil
$\Bl_{p_0} X \rightarrow \proj^1$ is generically an elliptic
fibration, with the proper transform $E_0$ of $p_0$ as zero-section.
An open subset of the moduli space of degree 1 Del Pezzo surfaces
correspond to fibrations where there are twelve singular fibers; in
this case the fibration is an elliptic fibration.  Let $\AAA$ be the
moduli space of such rational elliptic fibrations, along with a {\em
frame} for the base $\proj^1$.
A straight-forward dimension count shows that $\dim \AAA = 11$:
16 from choosing a pencil of plane cubics, +3 for a frame for the pencil,
-8 from $\Aut(\proj^2)$.

The structure of the group of sections of a rational elliptic
fibration is isomorphic to the $E_8$-lattice (via the Mordell-Weil height pairing).
There are 240 sections disjoint from $E_0$, corresponding to the 240
minimal vectors of the lattice.  Given any nine mutually disjoint
sections including $E_0$, those sections can be blown down, expressing
the fibration as the total space of plane cubics with distinct base
points (corresponding to the erstwhile sections).  There are $\#
W(E_8) = 2^{14} 3^5 5^2 7$ ways to choose the (ordered) octuple of disjoint
sections (not including than $E_0$).  Conversely, any pencil of cubics
with 12 nodal fibers has distinct base points, and the total space is
a rational elliptic fibration.  From this description, $\AAA$ is
smooth.  

Given a representation of a blow-up of a degree 1 Del Pezzo surface $X$ as the total space
of a pencil of cubics, we use the following convention for the homology of $X$:
$H$ is the hyperplane class of the plane, and $E_0$, \dots, $E_8$ are the base
points of the pencil (with $E_0$ corresponding to the exceptional divisor
of the blow-up of the Del Pezzo surface).  

We will make repeated use of the following fact about 2-torsion on a smooth elliptic fibration.

\tpoint{Lemma} \label{bungee} 
{\em Suppose $B$ is a curve, $b \in B$ is a closed point, and that $X
\rightarrow B$ is an elliptic fibrational with smooth total space,
such that the fiber $X_b$ above $b$ has multiplicative reduction.  Let
$\cL$ be a line bundle on $X$ of relative degree 2, and let $C$ be the
closure of the locus in $X$ of points $p$ such that in each smooth
fiber, $\oh(2p)$ is linearly equivalent to the restriction of $\cL$.
Then $C \rightarrow B$ is simply branched above $b$ (i.e.  the total
ramification index above $b$ is 1), and the branch point is the node
of $X_b$.}

\bpf
The 4-section $C$ meets $X_b$ at two points away from the node, and
hence with multiplicity two at the node.  Thus $C$ must be smooth at
the node, and $C \rightarrow B$ is branched with ramification index 1
there.  \epf

\bpoint{Weighted projective space}
\label{defB}
Consider the weighted projective space
$$
\proj  = \proj(  \underbrace{3, \dots, 3}_5, \underbrace{2, \dots, 2}_7)
$$
with projective coordinates $[a_0; \dots; a_4; b_0; \dots; b_6]$.  Let $\BBB$ be the 
open subscheme of $\proj$ where
the homogeneous degree 12 polynomial
$$
(a_0 x^4 + a_1 x^3 y + \dots + a_4 y^4)^3 + ( b_0 x^6 + b_1 x^5 y + \dots   + b_6 y^6)^2
$$ 
has distinct roots.  Then $\BBB$ is non-empty (as $[1,0,\dots,0,1] \in
\BBB$), and is clearly of dimension 11, nonsingular (as the only
singular points of $\proj$ are where $a_0 = \dots = a_4 = 0$ or $b_0 =
\dots = b_6 = 0$), and rational.

\bpoint{Branched covers of the projective line}
Recall that a {\em simply branched} cover of $\proj^1$ is defined to
be one where the total ramification index is 1 above each point of
$\proj^1$; the Hurwitz scheme is a moduli space for simply branched
connected covers.  The loci $\CCC$, $\DDD$, and $\EEE$ defined below
are locally closed subvarieties of the Hurwitz scheme.  By the
Riemann-Hurwitz formula, there are 12 branch points in each case,
giving a morphism (via the branch map) to the space of twelve points
in $\proj^1$.

\epoint{In the locus of genus 3 degree 4 covers of $\proj^1$, there is a codimension 1 set mapped to $\proj^1$ canonically}
Let $\CCC$ be the locus of smooth genus 3 degree 4 covers of $\proj^1$
branched over twelve distinct points, such that the pullback of
$\oh_{\proj^1}(1)$ is (isomophic to) the canonical bundle.  We call
this a {\em canonical pencil}.

\point \label{canorhyp}  Note that the source curve $C$ cannot be hyperelliptic.  Otherwise,
if $h: C \rightarrow \proj^1$ is the hyperelliptic map, then the
canonical bundle is (isomorphic to) $h^* \oh_{\proj^1}(2)$, and the
sections are pullbacks of sections of $h^* \oh_{\proj^1}(2)$.  Thus
any canonical pencil on $C$ factors through some 2-to-1 $j: \proj^1
\rightarrow \proj^1$; as $j$ is branched above 2 points, $h \circ j$
is not simply branched above those two points.

There are various equivalent formulations of the definition of $\CCC$.
For example, let $V$ be the (rank 3) Hodge bundle over $\cm_3$.  Let
$\grass_f(2,V)$ be the space of {\em framed} pencils in $\proj V$ over
$\cm_3$; $$
\dim \grass_f(2,V) = \dim \cm_3 + \dim G(2,3) + \dim \Aut \proj^1 = 11.
$$ Then $\CCC$ is an open substack (actually a scheme, as we saw
above) of $\grass_f(2,V)$ and rational (as $\cm_3$ is
rational, \cite{m3}), smooth, and of dimension 11.

An alternative formulation comes from the fact that if $C$ is not
hyperelliptic, then the canonical model of $C$ is a smooth quartic
curve in $\proj^2$.  Thus $\CCC$ is a quotient of an open subset of
$\proj H^0(\proj^2, \oh_{\proj^2}(4)) \times \proj^2$ by $\Aut \proj^2
= PGL(3)$ (with the additional data of a frame for the pencil): the
locus in $\proj H^0(\proj^2, \oh_{\proj^2}(4)) \times
\proj^2$ is the choice of a smooth quartic $Q \subset \proj^2$, and a
point $p$ not on any flex line or bitangent line of $Q$, nor on $Q$.
Then the corresponding cover of $\proj^1$ can be recovering by
projecting $Q$ from $p$.  In this guise, $\CCC$ is clearly smooth of
dimension 11, but rationality is not obvious.

\epoint{In the locus of genus 3 degree 4 covers of $\proj^1$, there is a codimension 1 set where the source curve is hyperelliptic}
Let $\DDD$ be the locus of smooth genus 3 degree 4 covers of $\proj^1$
branched over twelve distinct points, such that the source curve is
hyperelliptic.  By Section \ref{canorhyp}, the pencil is not a
canonical pencil (i.e. the pullback of $\oh_{\proj^1}(1)$ is not the
canonical bundle).  It is not hard to see that $\DDD$ is smooth of
dimension 11: the hyperelliptic locus is smooth of dimension 5; the
choice of non-canonical degree 4 line bundle gives 3 more dimensions;
such a line bundle has exactly 1 pencil by Riemann-Roch; and there are
3 dimensions of choice of frame for the pencil.

\epoint{In the locus of genus 4 degree 3 covers of $\proj^1$, there is a codimension 1 set where the pullback of $\oh_{\proj^1}(1)$ is
a theta-characteristic} In $\cm_4$, there is a divisor $\cm^1_4$
corresponding to curves $C$ with theta-characteristics with 2 sections
(vanishing theta-characteristic), smooth away from the hyperelliptic
locus $\ch_4$ (where it has 10 sheets corresponding to the Weierstrass
points).  Note that no (smooth) hyperelliptic genus 4 curve $C$ can
have a theta-characteristic inducing a base-point-free pencil, or
indeed any other degree 3 base-point-free pencil.  (Otherwise, let
$D_1$ be the divisor class of such a pencil, and let $D_2$ be the
hyperelliptic divisor class.  Define $$\phi: C \stackrel {(D_1,D_2)}
\longrightarrow \proj^1
\times \proj^1.$$  Then $\phi_*[C]$ is in class (2,3) on $\proj^1
\times \proj^1$, so $\phi$ cannot carry $C$ multiply onto its image.
Hence $\phi$ carries $C$ birationally onto its image.  But any curve
in class (2,3) has arithmetic genus 2, and the genus of $C$ is 4,
giving a contradiction.)

A point of $\cm^1_4 \setminus \ch_4$ corresponds to a curve $C$ whose 
canonical model is the intersection of a cone and a cubic
hypersurface in $\proj^3$; the pencil is given by the $\proj^1$ parametrizing the rulings of
the cone.
Let $\EEE$ be the open subset of such pencils, with framings, such
that the induced triple cover of $\proj^1$ is (simply) branched at 12
distinct points.  It is not hard to see that $\EEE$ is non-empty; in
any case, it will follow from Proposition \ref{EA}.  Clearly $\EEE$
is smooth of dimension 11, and as $\cm^1_4$ is rational (\cite{m14} p. 14),
$\EEE$ is as well.

\bpoint{Twelve points on the projective line:  polynomials that are a cube plus a square}  
All of the loci described above have natural morphisms to the moduli
space of twelve points on the projective line, $\Sym^{12} \proj^1
\cong \proj^{12}$, and all have the same image.  For concreteness, let
$\ZZZ$ be the locus of points in $\Sym^{12} \proj^1$ (considered as
homogeneous degree 12 polynomials in $x$ and $y$) corresponding to
polynomials that can be expressed as the sum of a square and a cube.
Then $\ZZZ$ is the image of $\BBB$ and hence locally closed.

\epoint{Other loci} 
\label{otherloci}
 There are various other descriptions of twelve points on $\proj^1$
that turn out to describe the same locus as those above.  As they are
either simple variants of the above, or are not well-defined moduli
problems, we will not dwell on them.
\begin{itemize}
\item In light of loci $\CCC$, $\DDD$, $\EEE$, a natural locus to consider
is that of hyperelliptic genus 5 curves branched over 12 points in the
locus $\ZZZ$.  A better description is to use the bijection between
simply branched double covers (see \cite{eehs} Theorem 6(1); when
$\proj^1$ is replaced by $\Q$, this is just part of the classical
method of solving the cubic), and \'{e}tale triple covers of double
covers (or equivalently a choice of 3-torsion of a double cover,
modulo $\{ \pm 1 \}$).  Under this ``elliptic-trigonal'' correspondence, there is a locus $\FFF$ corresponding to $\EEE$;
via the isomorphism of $\EEE$ with $\BBB$, it is essentially the
moduli space of representations of hyperelliptic genus 5 curves as
\begin{equation}
\label{genus5}
z^2 = f^2 + g^3
\end{equation}
where $f$ and $g$ are binary forms in $x$ and $y$ of degree 6 and 4
respectively.  (Note that by specifying the form of the equation
(\ref{genus5}) of the genus 5 curve, one automatically specifies an
\'{e}tale triple cover.)  The author is unaware of a geometric way of
describing this locus.
\item Project a smooth cubic surface in $\proj^3$ from a general point $p$, and let $B$ be the 
branch divisor in $\proj^2$.  ($B$ turns out to be a genus 4 curve
mapped canonically to $\proj^2$.)  Any conic in the plane meets the
sextic at twelve points, and hence describes a set of 12 points on $\proj^1$.  (This is related to $\EEE$ as follows: let
$C$ be the intersection of the cubic surface and the cone over the
conic with vertex $p$.)
\item In \cite{z2} Section 8, 
Zariski describes a locus of sextic plane curves $B'$, smooth except
for six cusps, and the six cusps lie on a conic.  Any conic meets this
sextic at twelve points.  (The sextic is the curve $B$ in
the previous example.  Zariski showed that the fundamental group of
the complement of $B'$ is $\zed/2 *
\zed/3$; the cubic surface can be recovered from the natural index 3
subgroup of this group.)
\end{itemize}

\bpoint{Relationship between moduli spaces}  
Surprisingly, it is possible to describe the relationship between
almost any pair of the spaces $\AAA$, $\BBB$, $\CCC$, $\DDD$, $\EEE$ (and through $\EEE$, $\ZZZ$)
described above.  Table 1 summarizes where the links are explained.

\begin{table}
\begin{center}
\begin{tabular}{|c|c|c|c|c|c|} \hline
       & $\AAA$ & $\BBB$ & $\CCC$ & $\DDD$ & $\EEE$     \\ \hline
$\AAA$ &   ---    &  \ref{ABCDE}   & \ref{CA}   &    \ref{DA}   &   \ref{EA}            \\ \hline 
$\BBB$ &   \ref{ABCDE}    &    ---   &  \ref{ganymede} & \ref{ganymede}       & \ref{ganymede}                \\ \hline
$\CCC$ &   \ref{CA}    &  \ref{ganymede}       &    ---   &        &  \ref{io}             \\ \hline
$\DDD$ &   \ref{DA}    &  \ref{ganymede}      &        &   ---   &   \ref{io}             \\ \hline
$\EEE$ &   \ref{EA}    &  \ref{ganymede}      &    \ref{io}   & \ref{io}   &    ---      \\ \hline
\end{tabular}
\end{center}
\caption{Sections where links are described}
\end{table}

We will repeatedly make use of the fact that if $\pi: X \rightarrow Y$
is a morphism of smooth varieties of the same dimension (over $\com$)
such that every point of $Y$ has $n$ pre-images, then $\pi$ is finite
\'{e}tale of degree $n$.

\tpoint{Theorem}  \label{ABCDE}
{\em There are natural isomorphisms $\AAA \cong \BBB \cong \EEE$, and
natural finite \'{e}tale morphisms $\CCC \rightarrow \AAA$ , $\DDD
\rightarrow \AAA$ of degree 120 and 135 respectively.}

The numbers 120 and 135 are related to the $E_8$ lattice; this
is made explicit in Section \ref{uranus}.

\bpf The isomorphsm $\AAA \cong \BBB$ follows from the classical theory of
Weierstrass models of rational elliptic fibrations (see for example
\cite{ms} Section 3, especially Theorem 1').  The rest of the
Theorem follows from Propositions \ref{CA} --- \ref{DA} below.   \epf

Let $M_1$ be the moduli space of degree 1 Del Pezzo surfaces.  If we
take quotients of $\AAA$, $\BBB$, and $\EEE$ by $PGL(2)$, we have
immediately:

\tpoint{Corollary}  
{\em $M_1$, $\BBB/PGL(2)$, and $\cm^1_4$ are birational (and the
isomorphic open sets are given quite explicitly).  As $\cm^1_4$ is
rational, $M_1$ and $\BBB/PGL(2)$ are as well.} 

Of course, this is well-travelled ground: $M_1$ is well-known to be
rational precisely because it is birational to $\cm^1_4$, see e.g.
\cite{m14} for example.  But Heckmann and Looijenga have
used the link to $\BBB$ to give a new proof of the rationality of 
$M_1$, or equivalently $\cm^1_4$ or $\BBB/PGL(2)$ (\cite{hl}).

We also get a degree 120 rational map from pencils in the Hodge bundle
over $\cm_3$ to $\cm^1_4$ (and an interpretation of this map
in terms of the $E_8$ lattice, see Section \ref{uranus}).

\bpoint{Ball quotients}
Deligne and Mostow (\cite{dm}) have shown that the moduli space of
unordered 12-tuples of points on a line (up to automorphisms of the
line $PGL(2)$) is a quotient of (a compactification of) the complex
9-ball by an (arithmetic) lattice in $PU(1,9)$.  Heckman and Looijenga
have shown (\cite{hl}) that the lift of this divisor in the locus of
twelve points is geodesically embedded in the 9-ball (and hence that
this locus is the quotient of a complex 8-ball by an arithmetic group.

\section{Relating branched covers to elliptic fibrations via 2-torsion}
\label{pluto}
The strategy we use to show a relationship between $\CCC$ (resp.
$\DDD$, $\EEE$) and $\AAA$ is as follows.  To construct an elliptic
fibration $X \rightarrow \proj^1$ (with zero-section $s$) from a
branched cover $C \rightarrow \proj^1$, we embed $C$ in a rational
surface $Y$ that is a $\proj^1$-bundle over $\proj^1$ (the Hirzebruch
surface $\eff_0$, $\eff_1$, or $\eff_2$), so that $C$ (possibly union
another curve) is a 4-section of the $\proj^1$-bundle.  We then
double-cover $Y$ branched over $C$ (and possibly another curve) to
construct $X$, a genus 1 fibration over $\proj^1$.  (In each case $X$
turns out to be rational.)  We then reconstruct the zero-section $s$
to express $X \rightarrow \proj^1$ as an elliptic fibration.

To construct a cover from an elliptic fibration $X \rightarrow
\proj^1$ (with zero-section $s$), we choose a line bundle on $X$ of
relative degree 2, and construct a 4-section that over each point over
$\proj^1$ consists of points $p$ such that $2p$ is linearly equivalent
to the restriction of the line bundle to the fiber.

\tpoint{Proposition}  \label{CA} 
{\em There is a natural degree 120 finite \'{e}tale morphism $\CCC \rightarrow \AAA$.}

We refer the interested reader also to \cite{do} Chapter IX.

\bpf  
We first describe the morphism $\CCC \rightarrow \AAA$.  For
simplicity of exposition, we define this morphism pointwise.  The data
of a point of $\CCC$ is a smooth quartic $Q$ in $\proj^2$ along with a
point $p \in \proj^2$ not on a bitangent line or flex line of $Q$, nor
on $Q$.  Let $Y = \Bl_p \proj^2$ (fibered over $\proj^1$ via the pencil
of lines through $p$, so $Y \cong \eff_1$); by abuse of notation we
consider $Q$ to be a subscheme of $Y$ (and a 4-section of the
fibration).

Let $X$ be the double cover of $Y$ branched over $Q$.  Such a double
cover exists, as $Q$ is divisible by 2 in $\Pic Y$.  (If the two
pre-images in $X$ of $p$ are blown down, by a classical result we have a
degree 2 Del Pezzo surface.  For a modern reference see \cite{cd}
Proposition 0.3.5.  The 56 (-1)-curves of the Del Pezzo surface as
the preimages of the 28 bitangents of $Q$.) Then $X$ is a rational
genus 1 fibration over $\proj^1$.  To make it an elliptic fibration,
we choose as zero-section either pre-image of the exceptional divisor
of $p$.  (The choice is irrelevant as they are interchanged under the
involution of the double cover $X \rightarrow Y$.)  The
nodal fibers of $X \rightarrow \proj^1$ correspond to the branch
points of $Q
\rightarrow \proj^1$, and indeed the nodes correspond to the
ramification points.

We next reverse the process.  Represent a point of $\AAA$ as a pencil
of cubics (with total space $X \rightarrow \proj^1$), with one of the
base points corresponding to $E_0$.  Then $E_1$ (corresponding to
another base point) is a section disjoint from $E_0$.  Let $E = E_2 +
\dots + E_8$ be the sum of exceptional divisors corresponding to the
remaining base points.  Let $Q'$ be the locus of points $p$ that, in
each fiber, satisfy $2p = E_0 + E_1$ (in the group law of the fiber).
Then $Q'$ is a smooth curve mapping 4-to-1 to $\proj^1$, simply
branched at 12 points (the nodes of the fibration, by Lemma \ref{bungee}); hence $Q'$ has
genus 3.  (It is not yet clear that $Q'$ is connected.) In $\Pic X$,
$$ Q' = aH + b(E_0+E_1) + cE $$ for some integers $a$, $b$, $c$.  As
$Q' \cap E_0 = Q' \cap E_1 = \emptyset$, $b=0$.  As $Q'$ is a 4-section
of the fibration, $$ 4 = Q' \cdot (3H-E_0-E_1 - E) = - Q' \cdot K_X =
3a + 7c.  $$ As $p_a(Q') = 3$, by adjunction, $$ 4 = Q' \cdot (Q' +
K_X) = (Q')^2 - 4, $$ so $8 = (Q')^2 = a^2-7c^2$.  Substituting
$3a+7c=4$, we find $(a-6)^2 = 0$, so $a=6$, $c=-2$.  By adjunction, $$
K_{Q'} = (3H-E)|_{Q'} = \oh_{\proj^1}(1)|_{Q'} $$ so the branched
cover $Q' \rightarrow \proj^1$ is indeed a canonical pencil.  Hence
$Q'$ is connected.  Observe that this construction commutes with the
morphism $\CCC \rightarrow \AAA$ described earlier.

The only choice involved was that of the section $E_1$ disjoint from
$E_0$; there are 240 such sections (Section \ref{ellfib}).  However,
they come in 120 pairs: if $(X, E_0, E_1) \rightarrow \proj^1$ is one
such pair, and $E_{-1}$ is the section corresponding to $-E_1$ in the
group law of the elliptic fibration, then there is an isomorphism
$$
\begin{array}{rcccl}
(X,E_0,E_1) & & \stackrel \sim {\longrightarrow} & & (X, E_{-1}, E_0) \\
& \searrow & & \swarrow \\
& & \proj^2
\end{array}
$$
(corresponding to translation by $E_{-1}$ in the group law of the elliptic fibration).  As the construction of $Q'$ 
above depended only on the data $(X, E_0 \cup E_1) \rightarrow \proj^1$, we see that for each point of $\AAA$ there
are 120 preimages in $\CCC$.  \epf

\tpoint{Proposition}   \label{EA} 
{\em There is a natural isomorphism $\EEE \stackrel \sim \rightarrow \AAA$.}

This is classical (see for example \cite{m14} p. 13 or \cite{cd} p. 36), but
we give a proof that highlights the importance of the genus 4 curve.

\bpf We first describe the morphism $\EEE \rightarrow \AAA$.  Consider a point of $\EEE$, or equivalently a
cone in $\proj^3$ transversely intersecting a cubic surface, along
with a frame for the $\proj^1$ parametrizing the rulings of the cone.
Let $Y$ be the blow-up of the cone at the vertex, isomorphic to the
Hirzebruch surface $\eff_2$.  Let $E$ be the exceptional divisor of
the blow-up.  Let $C$ be the intersection of the cone and the cubic,
pulled back to $Y$.  The induced morphism $C \rightarrow \proj^1$ is
precisely the triple cover in the definition of $\EEE$.  Let $X$ be the
double cover of $Y$ branched over $C+E$.  Such a cover exists as $C+E$
is even in $\Pic Y$.  (Reason: If $h$ is the pullback of
$\oh_{\proj^3}(1)$ to $Y \cong \eff_2$, and $f$ is a ruling of
$\eff_2$, then $C+E = (3h) + (h-2f) = 2 (2h-f)$.)  Then $X$ is a genus
1 fibration.  The choice of (preimage of) $E$ as 0-section makes $X$
an elliptic fibration.  Moreover, the blowdown of $X$ along $E$ is a
degree 1 Del Pezzo surface (classical; see \cite{cd} Proposition 0.3.6
for a modern reference).  Finally, the nodes of the fibers coincide
with the ramification points of $C \rightarrow \proj^1$, so $X
\rightarrow
\proj^1$ corresponds to a point of $\AAA$.

To reverse the process, represent a point of $\AAA$ as the total space
of a pencil of cubics $X \rightarrow \proj^1$, with one of the base
points corresponding to $E_0$.  Let $E = E_1 + \dots + E_8$ be the sum
of the exceptional divisors corresponding to the remaining base
points.  Let $C''$ be the closure of the locus of points $p$ that, in each fiber,
satisfy $2p = 2 E_0$.  Then by Lemma \ref{bungee} $C''$ is a smooth curve with $E_0$ as a
component; let $C' = C'' \setminus E_0$.  Thus $C'$ is a smooth curve
mapping 3-to-1 to $\proj^1$, simply branched at 12 points (the nodes
of the fibration); hence $C'$ has arithmetic genus 4 (although it is not yet
clear that $C'$ is connected).

In $\Pic X$, $C' = aH + bE + cE_0$.  As $C' \cap E_0 = \emptyset$,
$c=0$.  As $C'$ is a 3-section of the fibration, $$ 3 = C' \cdot
(3H-E_0-E) = - C' \cdot K_X = 3a+8b.  $$ As $p_a(C')=4$, by
adjunction, $$ 6 = (C')^2 + K_X \cdot C' = (C')^2 - 3, $$ so $a^2 -
8b^2=9$, from which $(a-9)^2=0$, so $a=9$, $b=-3$.  By adjunction,
$K_{C'} = 2(3H-E)|_{C'} = \oh_{\proj^1}(2)|_{C'}$, so the morphism $C'
\rightarrow \proj^1$ indeed corresponds to a theta-characteristic.
Hence $C'$ is connected.  Observe that this construction commutes with
the morphism $\EEE \rightarrow \AAA$ described earlier.
\epf

In the course of the proof, we recovered the (presumably classical)
result that the non-trivial 2-torsion points of fibers of a rational
elliptic fibration are transitively permuted by monodromy.

\tpoint{Proposition} \label{DA} 
{\em There is a natural degree 135 finite \'{e}tale morphism $\DDD
\rightarrow \AAA$.}

This proposition was known to Mumford (in the guise of the morphism
$\DDD \rightarrow \EEE$); see the discussion after Proposition
\ref{mercury}.

\bpf 
We begin by describing the ``reverse'' of the morphism 
$\DDD \rightarrow \AAA$, i.e., given a point of $\AAA$
and some choices,  how to construct a point of $\DDD$.

Consider a point of $\AAA$ represented as a pencil of cubics $X
\rightarrow \proj^1$ with section $E_0$ corresponding to one of the
base points (and $E_1$, \dots, $E_8$ correspond to the rest).  Then $H-E_0$ is a divisor on $X$ of relative degree 2.
Let $C'$ be the closure of the locus of points $p$ that, in each fiber, satisfy $2p =
H-E_0$ (in the group law of the fiber).  Then in $\Pic X$, $C' = aH +
b E_0 + cE$, where $E=E_1 + \dots + E_8$.  One can check (as in the
proofs of Propositions \ref{CA} and \ref{EA}) that $a=5$, $b=-3$,
$c=-1$, so the image of $C'$ in $\proj^2$ is a quintic with a triple
point (at the image $p_0$ of $E_0$).  Thus $C'$ is visibly
hyperelliptic (just project from $p_0$), so we have described a point
of $\DDD$.

(The value $b=-3$ corresponds to the enumerative fact that three times
in the pencil of cubics there is a flex at $p_0$; equivalently, there
are three fibers of the elliptic fibration where the divisor $H-3E_0$
is trivial.  We will use this fact in the discussion after Theorem
\ref{2tor}.)

By examining this construction, we will see how to construct the morphism
$\DDD \rightarrow \AAA$.  Consider the morphism
$$
\phi:  X \rightarrow \proj^1 \times \proj^1$$
given by the divisors $(H-E_0, 3H-(E_0+\dots + E_8))$;
$\phi$ maps $C$ to the first $\proj^1$ via the
hyperelliptic map.  Clearly $\phi$ is surjective of degree
2; it is branched over $C$.  The image
of $E_i$ ($1 \leq i \leq 8$) is a fiber of the first
projection; we see as
an added bonus that the points $E_i \cap C$  are the Weierstrass
points of $C$.

Thus any double-cover of $\proj^1 \times \proj^1$ branched over a
smooth divisor of type $(4,2)$ is a rational genus 1 fibration (with
the second projection as the morphism of the fibration, assuming that
the second morphism expresses the branch locus as a simply branched
cover of $\proj^1$), as this property is preserved under deformation.

We are now ready to describe the morphism $\DDD \rightarrow \AAA$.  Consider a
point of $\DDD$, or equivalently a hyperelliptic genus 3 curve $C$
(with hyperelliptic divisor class $D_1$), a degree 4 divisor class
$D_2$ (not the canonical divisor class), and a frame for the pencil of
sections of $D_2$.  Then there is a morphism $$
\phi:  C \stackrel { (D_1, D_2)} \longrightarrow \proj^1 \times \proj^1,
$$ determined up to automorphisms of the first $\proj^1$.  Let $Y$ be
the target $\proj^1 \times \proj^1$, and let $pr_i$ be the projection
to the $i$th factor.  Consider $Y$ as a fibration over (framed)
$\proj^1$ via the second projection $pr_2$.  Then $\phi$ is a closed
immersion (first note that $C$ is birational onto its image, then that
the arithmetic genus of the image is 3); for convenience, we identify
$C$ with its image.

Let $X$ be the double cover of $Y$ branched over the image of $C$; via
the fibration $pr_2: Y \rightarrow \proj^1$, $X$ is a genus 1 rational
fibration over $\proj^1$.  We now recover
a distinguished section $E_0$ to express $X$ as an elliptic fibration.

Let $F_1$, \dots, $F_8$ be the fibers of $\proj^1$ passing through the
hyperelliptic points of $C$.  The preimages of $F_i$ on $X$ are two
rational curves (meeting at a node, the hyperelliptic point of $C$);
call them $E_i$ and $E'_i$.  
By making these labellings, we have made $8! \times 2^7$ choices.  (We have
actually made $8! \times 2^8$ choices.  But exchanging $E_i$ with $E'_i$ for all $i$ gives
the same configuration up to the involution of the double cover $X \rightarrow Y$.)
Note that $E_i \cdot E_j = 0$ if $i \neq j$,
and 
$$
E_i^2 = E_i \cdot (E_i + E'_i - E'_i) = E_i \cdot ( pr_1^*(pt) - E'_i) = -1.
$$
Having chosen $E_1$, \dots, $E_8$, we can now determine $E_0$ uniquely.
The orthogonal complement of $\{ E_1, \dots, E_8 \}$ in $H^2(X,\zed)$ is $\zed H \oplus \zed E_0$; by
diagonalizing (to get $\left[ \begin{matrix} 1 & 0 \\ 0 & -1 \end{matrix} \right]$) we can find $E_0$.
(More precisely, diagonalize to determine the class of $E_0$ in
$H^2(X,\zed)$ up to sign.  Then note that $E_0$ is effective to
determine the class of $E_0$ in $H^2(X,\zed)$.  $E_0$ is the only
effective representative in its class.)

In short, given a rational elliptic fibration with choice of sections $E_1$, 
\dots, $E_8$ (there are $\# W(E_8)$ possible choices by Section \ref{ellfib}), we
recover a point of $D$ along with one of $8! \times 2^7$ choices.  Hence
each point of $\DDD$ has
$$
\# W(E_8) / (8!  \times 2^7) = 135
$$
preimages.
\epf

\section{Recillas' trigonal construction, and theta-characteristics of the genus 4 curve}
\label{uranus}
We recall Recillas' beautiful construction (\cite{r}) giving a
bijection between connected simply branched quadruple covers of
$\proj^1$ to connected unramified double covers of connected simply
branched triple covers of $\proj^1$.  The triple cover is simply
branched over the same points as the quadruple cover.  Furthermore,
the Jacobian of $Q$ is isomorphic to the Prym of $S \rightarrow C$.
See \cite{donagi} for an excellent exposition (the Prym result is a
variant of Donagi's tetragonal construction).

For simplicity, we give a description of the bijection over the complex numbers, although the
construction works over any algebraically closed field of characteristic at least 5. 
Given the data of a simply branched quadruple cover $Q \rightarrow \proj^1$, one recovers $S \rightarrow C \rightarrow \proj^1$ (where
$S \rightarrow C$ is an unramified double cover, and $C \rightarrow \proj^1$ is a simply branched triple cover) as follows.
After choosing a base point and making branch cuts to the branch points of $Q$, label the sheets of $Q$ $a, b, c, d$  
(and remember how monodromy around the branch points of $Q \rightarrow \proj^1$ permutes the sheets).  Then $C \rightarrow \proj^1$
is obtained by considering three sheets labelled 
$\{ \{ a, b \}, \{c, d \} \}$,
$\{ \{ a, c \}, \{b,  d \} \}$, 
$\{ \{ a, d \}, \{b, c \} \}$ (each is a set of sets; for simplicity call these $ab+cd$, $ac+bd$, and $ad+bc$), where
the action of the monodromy group on the three sheets is as induced by the action on the four sheets of $Q$.
Also, $S \rightarrow \proj^1$ is obtained by considering six sheets labelled
$\{ a, b \}$, $\{c, d \}$,
$ \{ a, c \}$, $\{ b,  d  \}$, 
$ \{ a, d \}$, $\{ b, c \}$; the morphism $S \rightarrow C$ is as one would expect.
One can quickly check that $C \rightarrow \proj^1$ is simply branched and that $S \rightarrow C$ is \'{e}tale.
This construction is summarized pictorially in Figure 1.

\begin{figure}
\begin{center}

	   \setlength{\unitlength}{.1\baseunit}
	    \begingroup\makeatletter\ifx\SetFigFont\undefined%
\gdef\SetFigFont#1#2#3#4#5{%
  \reset@font\fontsize{#1}{#2pt}%
  \fontfamily{#3}\fontseries{#4}\fontshape{#5}%
  \selectfont}%
\fi\endgroup%
{\renewcommand{\dashlinestretch}{30}
\begin{picture}(10867,4614)(0,-10)
\path(3612,3687)(12,3687)
\path(6612,4287)(10212,4287)
\path(6612,4587)(10212,4587)
\path(6612,3987)(10212,3987)
\path(6612,3687)(10212,3687)
\path(6612,3387)(10212,3387)
\path(6612,3087)(10212,3087)
\path(6612,2187)(10212,2187)
\path(6612,1887)(10212,1887)
\path(6612,1587)(10212,1587)
\path(3612,3387)(12,3387)
\path(3612,3087)(12,3087)
\path(3612,2787)(12,2787)
\path(7212,12)(3612,12)
\path(1812,2562)(5112,312)
\blacken\path(4995.953,354.814)(5112.000,312.000)(5029.753,404.387)(5042.597,359.320)(4995.953,354.814)
\path(8412,1362)(5712,312)
\blacken\path(5812.967,383.454)(5712.000,312.000)(5834.714,327.533)(5790.288,342.445)(5812.967,383.454)
\path(8412,2937)(8412,2337)
\blacken\path(8382.000,2457.000)(8412.000,2337.000)(8442.000,2457.000)(8412.000,2421.000)(8382.000,2457.000)
\put(3912,3012){\makebox(0,0)[lb]{\smash{{{\SetFigFont{8}{9.6}{\rmdefault}{\mddefault}{\updefault}$c$}}}}}
\put(3912,2712){\makebox(0,0)[lb]{\smash{{{\SetFigFont{8}{9.6}{\rmdefault}{\mddefault}{\updefault}$d$}}}}}
\put(10437,4512){\makebox(0,0)[lb]{\smash{{{\SetFigFont{8}{9.6}{\rmdefault}{\mddefault}{\updefault}$ab$}}}}}
\put(10437,4212){\makebox(0,0)[lb]{\smash{{{\SetFigFont{8}{9.6}{\rmdefault}{\mddefault}{\updefault}$cd$}}}}}
\put(10437,3912){\makebox(0,0)[lb]{\smash{{{\SetFigFont{8}{9.6}{\rmdefault}{\mddefault}{\updefault}$ac$}}}}}
\put(10437,3612){\makebox(0,0)[lb]{\smash{{{\SetFigFont{8}{9.6}{\rmdefault}{\mddefault}{\updefault}$bd$}}}}}
\put(10437,3312){\makebox(0,0)[lb]{\smash{{{\SetFigFont{8}{9.6}{\rmdefault}{\mddefault}{\updefault}$ad$}}}}}
\put(10437,3012){\makebox(0,0)[lb]{\smash{{{\SetFigFont{8}{9.6}{\rmdefault}{\mddefault}{\updefault}$bc$}}}}}
\put(10437,1812){\makebox(0,0)[lb]{\smash{{{\SetFigFont{8}{9.6}{\rmdefault}{\mddefault}{\updefault}$ac+bd$}}}}}
\put(10437,1512){\makebox(0,0)[lb]{\smash{{{\SetFigFont{8}{9.6}{\rmdefault}{\mddefault}{\updefault}$ad+bc$}}}}}
\put(3912,3612){\makebox(0,0)[lb]{\smash{{{\SetFigFont{8}{9.6}{\rmdefault}{\mddefault}{\updefault}$a$}}}}}
\put(3912,3312){\makebox(0,0)[lb]{\smash{{{\SetFigFont{8}{9.6}{\rmdefault}{\mddefault}{\updefault}$b$}}}}}
\put(10437,2112){\makebox(0,0)[lb]{\smash{{{\SetFigFont{8}{9.6}{\rmdefault}{\mddefault}{\updefault}$ab+cd$}}}}}
\end{picture}
} 
\end{center}
\label{Recillas1}
\caption{One direction of the Recillas correspondence}
\end{figure}
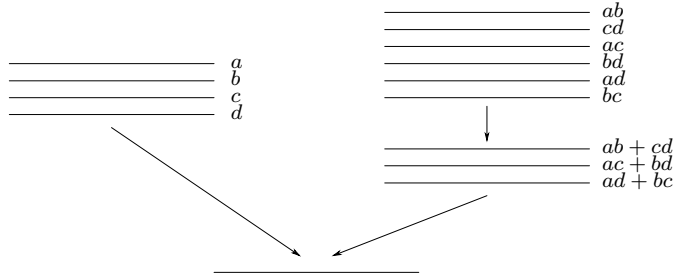

The opposite direction is more elaborate.  Let the three sheets of $C$
be called $x$, $y$, and $z$, and let the six sheets of $S$ be called
$x'$, $x''$, $y'$, $y''$, $z'$, $z''$ (with $x'$ and $x''$ the
preimages of of $x$, etc.).  Then the monodromy group of $S
\rightarrow \proj^1$ is a subgroup of the symmetries of an octahedron:
consider an octahedron with vertices $x'$, \dots, $z''$, with $x'$
opposite $x''$, $y'$ opposite $y''$, and $z'$ opposite $z''$.  Then
the octahedron has four pairs of opposite faces; this gives the four sheets
of $Q$ over $\proj^1$.  This construction is summarized pictorially in Figure 2.

\begin{figure}
\begin{center}

	   \setlength{\unitlength}{.1\baseunit}
	    \begingroup\makeatletter\ifx\SetFigFont\undefined%
\gdef\SetFigFont#1#2#3#4#5{%
  \reset@font\fontsize{#1}{#2pt}%
  \fontfamily{#3}\fontseries{#4}\fontshape{#5}%
  \selectfont}%
\fi\endgroup%
{\renewcommand{\dashlinestretch}{30}
\begin{picture}(11388,4614)(0,-10)
\path(12,4287)(3612,4287)
\path(12,4587)(3612,4587)
\path(12,3987)(3612,3987)
\path(12,3687)(3612,3687)
\path(12,3387)(3612,3387)
\path(12,3087)(3612,3087)
\path(12,2187)(3612,2187)
\path(12,1887)(3612,1887)
\path(12,1587)(3612,1587)
\path(1812,2937)(1812,2337)
\blacken\path(1782.000,2457.000)(1812.000,2337.000)(1842.000,2457.000)(1812.000,2421.000)(1782.000,2457.000)
\path(7212,12)(3612,12)
\path(1812,1362)(5112,312)
\blacken\path(4988.553,319.797)(5112.000,312.000)(5006.745,376.972)(5031.954,337.469)(4988.553,319.797)
\path(8412,2562)(5712,312)
\blacken\path(5784.981,411.869)(5712.000,312.000)(5823.392,365.775)(5776.531,365.775)(5784.981,411.869)
\path(10212,4137)(6612,4137)
\path(10212,3687)(6612,3687)
\path(10212,3237)(6612,3237)
\path(10212,2787)(6612,2787)
\put(3912,4512){\makebox(0,0)[lb]{\smash{{{\SetFigFont{8}{9.6}{\rmdefault}{\mddefault}{\updefault}$x'$}}}}}
\put(3912,4212){\makebox(0,0)[lb]{\smash{{{\SetFigFont{8}{9.6}{\rmdefault}{\mddefault}{\updefault}$x''$}}}}}
\put(3912,3912){\makebox(0,0)[lb]{\smash{{{\SetFigFont{8}{9.6}{\rmdefault}{\mddefault}{\updefault}$y'$}}}}}
\put(3912,3612){\makebox(0,0)[lb]{\smash{{{\SetFigFont{8}{9.6}{\rmdefault}{\mddefault}{\updefault}$y''$}}}}}
\put(3912,3312){\makebox(0,0)[lb]{\smash{{{\SetFigFont{8}{9.6}{\rmdefault}{\mddefault}{\updefault}$z'$}}}}}
\put(3912,3012){\makebox(0,0)[lb]{\smash{{{\SetFigFont{8}{9.6}{\rmdefault}{\mddefault}{\updefault}$z''$}}}}}
\put(3912,2112){\makebox(0,0)[lb]{\smash{{{\SetFigFont{8}{9.6}{\rmdefault}{\mddefault}{\updefault}$x$}}}}}
\put(3912,1812){\makebox(0,0)[lb]{\smash{{{\SetFigFont{8}{9.6}{\rmdefault}{\mddefault}{\updefault}$y$}}}}}
\put(3912,1512){\makebox(0,0)[lb]{\smash{{{\SetFigFont{8}{9.6}{\rmdefault}{\mddefault}{\updefault}$z$}}}}}
\put(10512,4062){\makebox(0,0)[lb]{\smash{{{\SetFigFont{8}{9.6}{\rmdefault}{\mddefault}{\updefault}$x'y'z'+x''y''z''$}}}}}
\put(10512,3612){\makebox(0,0)[lb]{\smash{{{\SetFigFont{8}{9.6}{\rmdefault}{\mddefault}{\updefault}$x''y'z'+x'y''z''$}}}}}
\put(10512,3162){\makebox(0,0)[lb]{\smash{{{\SetFigFont{8}{9.6}{\rmdefault}{\mddefault}{\updefault}$x'y''z'+x''y'z''$}}}}}
\put(10512,2712){\makebox(0,0)[lb]{\smash{{{\SetFigFont{8}{9.6}{\rmdefault}{\mddefault}{\updefault}$x'y'z''+x''y''z'$}}}}}
\end{picture}
} 
\end{center}
\label{Recillas2}
\caption{The other direction of the Recillas correspondence}
\end{figure}
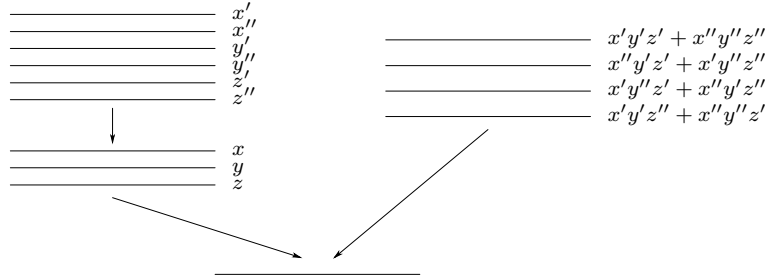

\tpoint{Proposition}\label{io}
{\em The morphisms $\CCC \rightarrow \EEE$ and $\DDD \rightarrow \EEE$
(obtained by composing the morphisms and isomorphisms of Propositions
\ref{CA}, \ref{EA}, and \ref{DA}) are obtained by Recillas' trigonal construction.}

\bpf
Consider a point of $\CCC$, corresponding to a quadruple (canonical)
cover $Q \rightarrow
\proj^1$, and any point $q$ of $\proj^1$ that is not a branch point of
$Q \rightarrow \proj^1$.  Following the proof of Proposition \ref{CA},
embed $Q$ in a rational elliptic fibration, where the points of $Q$
above $q$ are distinct points $a$, $b$, $c$, $d$ such that
$2a=2b=2c=2d$ in the group law of the fiber.  Then to each unordered
pair of elements of $\{a, b, c, d \}$ we can associated their
difference, which is a non-zero two-torsion point of the fiber (the
order of the pair doesn't matter).  Also, the complementary pair is
associated to the same two-torsion point.  Hence we have identified
the set $\{ ab+cd, ac+bd, ad+bc \}$ (which appears in Recillas'
construction) with the set of non-zero two-torsion points of the fiber
(which agrees with the description of $\EEE$ given in Proposition
\ref{EA}).

The proof for $\DDD$ is identical.
\epf

The morphism can be related to the beautiful geometry of the genus 4
curve in the description of $\EEE$.  Consider a given point of $\EEE$,
corresponding to a curve $C$ with vanishing theta-characteristic
$\theta$ inducing the base-point-free pencil $C \rightarrow \proj^1$.
As $C$ has genus 4, $C$
has exactly $\# H^1(C, \zed/2)-1 = 255$ connected \'{e}tale double
covers, so all such are accounted for by the tetragonal covers of the
form $\CCC$ and $\DDD$.

Now $C$ has 120 odd theta-characteristics and 136 even ones (including
$\theta$).  Translating by $\theta$, we have 120 ``odd'' 2-torsion
points of $\Pic C$ (all non-zero), and 135 non-zero ``even'' 2-torsion
points.  The non-zero 2-torsion points parametrize double covers of
$C$.  Globalizing this construction over $\EEE$, we have an \'{e}tale
degree 255 morphism $\EEE' \rightarrow \EEE$; $\EEE'$ is disconnected,
and splits into two pieces (not a priori connected) $\EEE'_{\odd}$
and $\EEE'_{\even}$, of degree 120 and 135 over $\EEE$ respectively.
Hence we must have $\EEE'_{\odd} = \CCC$ and $\EEE'_{\even} = \DDD$.  As a consequence
we have:

\tpoint{Proposition}  
\label{mercury} 
{\em
\begin{itemize}
\item[(a)]  The monodromy group of odd (resp. even) non-zero 2-torsion of the Picard group of the universal 
curves over $\cm^1_4 \setminus \ch_4$ is full, i.e.
$\# W(E_8) / \{ \pm 1 \}$.  The same statement is true with non-zero 2-torsion of Picard replaced by
the non-vanishing theta-characteristics.
\item[(b)]  The Jacobian of any genus 3 curve $C_\CCC$ is the Prym of a rational 2-parameter family of genus 4
curves with vanishing theta-characteristic (parametrized by an open subset of $\proj H^0(C_\CCC, K_{C_\CCC})^\vee$).
\item[(c)]  The Jacobian of a hyperelliptic curve $C_\DDD$ is the Prym of a 3-parameter family of genus 4
curves with vanishing theta-characteristic (parametrized by an open subset of $\Pic^4 C_{\DDD}$).
\end{itemize}
}

\bpf 
All that remains to be proved is the monodromy statement on
theta-characteristics (or equivalently 2-torsion of Picard) on the
universal curve over $\cm^1_4 \setminus \ch_4$.  But if $L$ is the universal
$E_8$-Mordell-Weil-lattice over $\AAA$, then the monodromy group on $L$ is full
(i.e. $W(E_8)$), and $L / L(2)$ is now identified with the 2-torsion
of the universal curve.
\epf

In light of the above discussion, we see that in Proposition \ref{DA} 
we have recovered a result of Mumford.  (We are grateful to R. Smith
for pointing this out.)  Mumford has shown (Theorem (c) of \cite{pv1}
p. 344): If $\tC \rightarrow C$ is an \'{e}tale double cover of a
non-hyperelliptic genus 4 curve $C$ (corresponding to $L \in
\Pic(C)[2]$, with Prym variety $P$ (with theta divisor $\Xi$), then
$(P,\Xi)$ is a (genus 3) Jacobian, and $\Xi$ is singular iff $(P,\Xi)$
is a hyperelliptic Jacobian iff $C$ has a vanishing
theta-characteristic $\theta$ such that $\theta \otimes L$ is an even
theta-characteristic.  (His method of proof is different.)

\bpoint{Explicit construction of 2-torsion of a genus 4 curve with vanishing theta-characteristic, and
relation to $\CCC$ and $\DDD$}
Given the explicitness of the above constructions, one should expect
to see all of the theta-characteristics and 2-torsion of a genus 4
curve in $\cm^1_4 \setminus \ch_4$ in a particularly straightforward way.  For concreteness,
fix such a curve $C$ (with theta-characteristic $\theta$ inducing
pencil $\pi: C \rightarrow \proj^1$) and embed it in the total space
$X$ of a pencil of cubics as described in the proof of Proposition
\ref{EA}.  Then $X$ is the blow-up of $\proj^2$ (with hyperplane class
$H$) at 9 points, with exceptional divisors $E_0$, \dots, $E_8$, and
$C$ is in class $9H-3(E_1  + \dots + E_8)$.  $X$ is an elliptic fibration
(with zero-section $E_0$), and $C$ is the non-zero 2-torsion of the fibration.
Note that there is a natural involution $\iota$ of $X \rightarrow \proj^1$ preserving $C$, which is
the inverse in the group law of the fibration.

\tpoint{Theorem}
\label{2tor}
{\em
The restriction map
\begin{equation}
\si:  \Pic X \rightarrow C
\end{equation}
surjects onto $\zed \theta \oplus \Pic C[2]$.}

We can describe $\sigma$ explicitly.  Note that $E_0 \cap C =
\emptyset$, so $\si(E_0) = 0$.  From the proof of Proposition \ref{EA},
$\si(3H-(E_1+ \dots + E_8)) = \theta$.

The genus 3 curve $C_\CCC$ in $X$ described in the proof of Proposition \ref{CA} lies
in class $6H - 2 (E_2 + \dots + E_8)$, so $C \cdot C_\CCC = 12$.  But
$C$ and $C_\CCC$ both pass through the 12 nodes of the elliptic fibration, at the points where both $\pi: C \rightarrow
\proj^1$ and $C_\CCC \rightarrow \proj^1$ ramify.  Hence $\si(C_\CCC)$ is the ramification divisor of $\pi$.
By the Riemann-Hurwitz formula,
$
K_C= \pi^* K_{\proj^1} + \si(C_\CCC)$, so 
$$
\si( 6H - 2 (E_1 + \dots + E_8)) = \si(-6H+2 (E_0 + \dots + E_8) + 6H - 2 (E_2 + \dots + E_8)
$$
$$
\Longrightarrow
\si( 6H - 2 (E_1 + \dots + E_8) - 2E_1) = 0.
$$ Hence $\si(3H - 2 E_1 - (E_2 + \dots + E_8)) \in \Pic C[2]$; this
corresponds to the double cover associated to the tetragonal
cover $C_\CCC \rightarrow \proj^1$.  (It would be interesting
to understand this correspondence more explicitly.)
Also $\theta = \si(3H - (E_1 + \dots +
E_8))$, so $\si(E_1)$ is a theta-characteristic, and by symmetry
$\si(E)$ is a theta-characteristic for any of the 240
exceptional curves on $X$ not meeting $E_0$ (i.e. the exceptional
curves on the degree 1 Del Pezzo surface that is $X$ blown down along
$E_0$).    As $\si(E) = \si( \iota (E))$ (where $\iota$ is the
involution described just before the statement of Theorem \ref{2tor}),
they come in pairs; they are the 120 odd theta-characteristics.

The hyperelliptic genus 3 curve $C_\DDD$ described in the proof of Proposition
\ref{DA} lies in class $5H-3E_0 - (E_1 + \dots + E_8)$, so $C_\DDD \cdot C = 21$.
Note that $C_\DDD$ and $C$ meet at the three non-zero 2-torsion points
in the three fibers where $H=3E_0$ (restricted to the fiber).  
(These three fibers were mentioned in the
proof of Proposition \ref{DA}.)  This is because in such a fiber, $H-E_0 =
2E_0$, and the restriction of $C_\DDD$ (resp. $C$) to the fiber are
those points $p$ such that $2p = H-E_0$ (resp.  $2p= 2E_0$ and $p \neq
0$) in the group law of the fiber.  Thus $C_\DDD$ meets $C$
transversely at these 9 points, and at the 12 nodes of the fibration, so
\begin{equation}
\label{star}
\si(C_\DDD) = \si(C_\CCC) + 3 \si( 3H - (E_0 + \dots + E_8))
\end{equation}
From the previous paragraph,
\begin{eqnarray*}
\si(C_\CCC) &=& \si \left(  (6H- 2 (E_2 + \dots + E_8) ) + ( 6H - 4E_1 - 2(E_2 + \dots + E_8)) \right) \\
&=& \si(12H-4 (E_1 + \dots + E_8))
\end{eqnarray*}
Substituting this into (\ref{star}):
$$
\si(16H - 6(E_1 + \dots + E_8)) = 0.
$$
Thus 
$\si(8H - 3(E_1 + \dots + E_8)) \in \Pic C[2]$, and
this corresponds to the double cover associated to the tetragonal curve $C_\DDD \rightarrow \proj^1$.

{\noindent {\em Proof of Theorem \ref{2tor}.}} $\si(E_0) = 0$,
$\si(3H-(E_1 + \dots + E_8))$ and $\si(E_i)$ ($i>0$) are
theta-characteristics, and $\si(8H-3(E_1+\dots + E_8)) \in \Pic C[2]$;
hence the image of $\sigma$ is contained in $\zed \theta \oplus \Pic
C[2]$.  On the other hand, as all theta-characteristics are in the
image of $\sigma$, the image of $\sigma$ contains $\zed \theta \oplus
\Pic C[2]$.  \epf

\bpoint{Summary:  A unified picture of $\AAA$, $\CCC$, $\DDD$, $\EEE$}
The results of this section can be summarized in a single picture
essentially due to Dolgachev and Ortland (\cite{do} Section VII.5).
Above the moduli space $\AAA$ of framed degree 1 Del Pezzo surfaces
$X$, there is a cover corresponding to the Mordell-Weil lattice modulo
even vectors.  This is the same cover corresponding to 2-torsion in
the Picard group of the genus 4 curve parametrized by $\EEE$.  The
author is grateful to N. Elkies for pointing out that the Weil-pairing
on 2-torsion corresponds in the second incarnation to half of the
Mordell-Weil pairing in the first incarnation.

This degree 256 cover splits into $\CCC$, $\DDD$, and the
zero-section.

\section{Explicit description of the morphisms $\CCC \rightarrow \BBB$, $\DDD \rightarrow \BBB$, $\EEE \rightarrow \BBB$ via discriminants}
\label{ganymede}
(This section is classical, although we hope the presentation is of interest.)
Although the morphisms $\CCC \rightarrow \BBB$, $\DDD \rightarrow \BBB$, $\EEE \rightarrow \BBB$ are now clear via
elliptic fibrations (Section \ref{pluto}), it is enlightening to see them in terms of the elementary algebra of discriminants.
Let $f(a,b)$ be the binary quartic with indeterminate coefficients
\begin{equation}
\label{comet}
f(a,b) = p_0 a^4 + p_1 a^3 b + p_2 a^2 b^2 + p_3 a b^3 + p_4 b^4.
\end{equation}
Then the discriminant of $f$ is $\Delta f = 4 u_2^3 + 27 u_3^2$ where
\begin{eqnarray*}
u_2 &=& p_1 p_3 - 4 p_0 p_4 - p_2^2/3, \\
u_3 &=& p_1^2 p_4 + p_0 p_3^2 - 8 p_0 p_2 p_4/3 -  p_1 p_2 p_3/3 + 2 p_2^3/27.
\end{eqnarray*}
If the $p_i$ are binary forms of degree 2 in $x$ and $y$, $\Delta f$ is of degree 12, and
(as it is a square plus a cube) corresponds to a point in $\BBB$.  Then (\ref{comet}) describes
a hyperelliptic curve mapping 4-to-1 onto $\proj^1$ (and 
specifically a (4,2)-class 
on $\proj^1 \times \proj^1$, which appeared in the proof of Proposition \ref{DA}).

If the $p_i$ are binary forms of degree $i$ in $x$ and $y$, then
$\Delta f$ is again of degree 12, and again corresponds to a point of
$\BBB$; then (\ref{comet}) describes a plane quartic with the data of
an additional point in the plane (more precisely a curve on the
Hirzebruch surface $\eff_1$), and we recover the construction of the
proof of Proposition \ref{CA}.

The multiplicities of 120 and 135 in the above cases are not clear.

The same idea works with cubics.  Let
\begin{equation}
\label{comet2}
f(a,b) = q_0 a^3 + q_1 a^2 b + q_2 a b^2 + q_3  b^3.
\end{equation}
Then the discriminant of $f$ is $\Delta f = (4 u_2^3 + 27 u_3^2)/q_0^2$ where
\begin{eqnarray*}
u_2 &=& q_0 q_2 -  q_1^2/3 , \\
u_3 &=& q_0^2 q_3 - q_0 q_1 q_2/3 +  2 q_1^3 /27.
\end{eqnarray*}
If the $q_i$ are binary forms of degree $2i$ in $x$ and $y$, and $q_0=1$, then
$\Delta f$ has degree 12, hence corresponds to a point of 
$\BBB$, and (\ref{comet}) describes the restriction of a cubic hypersurface
to a quadric cone in $\proj^3$, so we recover the construction of the
proof of Proposition \ref{EA}.  

\section{Twelve points on the projective line:  the locus $\ZZZ$}
\label{saturn}

\tpoint{Theorem} \label{bz} 
{\em \begin{enumerate}
\item[(a)]  The morphism $\pi:  \BBB \rightarrow \ZZZ$ is unramified and proper.
\item[(b)]  $\pi$ is birational, but not an isomorphism.
\item[(c)]  $\pi$ is the normalization of $\ZZZ$.
\end{enumerate}
}

(Clearly (c) implies (b); however, we use (b) to prove (c).)

Hence $\ZZZ$ is rational of dimension 11.  Also, this result shows
that the title of the paper is somewhat misleading: the locus $\ZZZ$
(of twelve points on the projective line) is the ``wrong'' moduli
space to study, and the ``right'' spaces are $\AAA$, $\BBB$, etc.

\bpf
(a) We need only show that $\pi: \BBB \rightarrow \Sym^{12} \proj^1$
is unramified and proper.  Properness is immediate: if $\proj$ is the
weighted projective space from the definition of $\BBB$, then $\proj
\rightarrow \Sym^{12} \proj^1$ is proper, and $\BBB$ was defined as the subset
of $\proj$ disjoint from the preimage of the discriminant locus $\De$,
so $\BBB \rightarrow \Sym^{12} \proj^1 \setminus \De$ is proper. 

To check that $\pi$ is unramified, note that by the
``hyperelliptic-trigonal'' correspondence of Section \ref{otherloci},
$\BBB$ can be associated with a locally closed subscheme of the
\'{e}tale cover $X \rightarrow \Sym^{12} \proj^1 \setminus \De$, where
$X$ is the non-zero three-torsion (modulo $\{ \pm 1 \}$) of the
hyperelliptic curve branched at those 12 points.

(One can also check that $\pi$ is unramified directly,
by explicitly describing a general point and tangent vector of 
$\BBB$, and computing the image in $\Sym^{12} \proj^1$.) 

(b)
Assume that $\BBB \rightarrow \ZZZ$ has degree greater than one.  Then as $\dim \BBB = 11$, the parameter space $P$  of homogeneous
polynomials $f_1$, $g_1$, $f_2$, $g_2$ (in two variables) such that
\begin{enumerate}
\item[(i)] $f_1^3 + g_1^2 = f_2^3 + g_2^2$,
\item[(ii)] $\deg f_i = 4$, $\deg g_i = 6$,
\item[(iii)] $g_1^2 \neq g_2^2$ (or equivalently $f_1^3 \neq f_2^3$)
\end{enumerate}
is of dimension at least 12.  But there are 10 dimensions of choices
of coefficients of $f_1$ and $f_2$.  Once $f_1$ and $f_2$ are given,
there is a one-parameter family of choices of $g_1$ and $g_2$ (through
$(g_2+g_1)(g_2-g_1) = f_1^3 - f_2^3 \neq 0$, as the
roots of $f_1^3-f_2^3$ must be split between $(g_2+g_1)$ and $(g_2-g_1)$; the one
parameter comes from the choice of leading coefficient of $(g_2+g_1)$).
Hence $\dim P = 11$, giving a contradiction.  Thus $\BBB$ is
birational.

As $(x^6)^2 + (y^4)^3 = (y^6)^2 + (x^4)^3$ 
has 12 distinct roots, $\BBB \rightarrow \ZZZ$ isn't injective and thus isn't an isomorphism.

(c) From (a), $\pi$ is proper and (as $\BBB$ is smooth and $\pi$ is unramified) quasifinite,
so $\pi$ is finite.  As $\pi$ is birational and $\BBB$ is normal,
$\pi$ is the normalization of $\ZZZ$.
\epf

As a side benefit, we see from the proof that at any point of $\ZZZ$, the branches are smooth.

Note that the analogous argument works for $\Sym^n \proj^1$ when
$n=6k$, $k \geq 2$ (i.e. a general homogeneous polynomial of degree
$n$ that is expressible as the sum of a cube and a square, is so
expressible in only one way), but the proof of (b) (and hence (c))
breaks down when $k=1$.  Indeed, the general sextic can be written as
a cube plus a square in 40 essentially different ways (\cite{e}
Theorem 3 i), a result of Clebsch).  Also, the proof that $\pi$ is
unramified in (a) must be done explicitly.

\tpoint{Theorem}
{\em The degree of $\ZZZ$ in $\proj^{12}$ is 3762.}
\label{degZ} 

\bpf
In \cite{quartic} Section 9.1, it was shown that the number of genus 3 canonical covers of $\proj^1$ (i.e.
points of $\CCC$) with
11 fixed branch points is $3762 \times 120$.  The result then follows Theorems \ref{ABCDE} and
\ref{bz}.  

Alternatively, Zariski computed degree $\deg \ZZZ = 3762$ via the locus $\FFF$ described in
Section \ref{otherloci}; see Section \ref{classicz}.
\epf

It would be interesting to derive this degree more directly.

\epoint{Remark}  W. Lang has proved that the degree of the locus $\AAA$ in 
characteristic 2 is 1 --- remarkably, the condition for twelve points
to appear as the discriminant locus of a rational elliptic fibration
is linear (\cite{bill}).

\section{Classical results}
\label{cr}

Some of the above links have been described classically.

\epoint{Zeuthen}  In \cite{zeuthen} (p. XXII), Zeuthen solves 
an enumerative problem that, in modern language, translates to: given
11 points on a line, how many canonical covers are there branched at
those 11 points?  He gives the correct answer ($451440 = 120 \times
3762$), but it is unclear how he obtained this.  (This fact falls out
as a side benefit of the proof of \cite{quartic} Theorem 8.1.)  More
precisely, he asks a different question, to which he gives an
incorrect answer, without throwing off his calculation of the
characteristic number of quartic curves.  The reason for his error, in
modern terms, is that he did not suspect that the degree of $\CCC
\rightarrow \ZZZ$ was 120.  This is discussed in \cite{quartic} Section
9.1.

\epoint{Zariski}
\label{classicz}
In \cite{z1}, Zariski calculates the answer to a similar question:
he defines the locus $\Gamma_{12}$ in $\Sym^{12} \proj^1$ as the image of
$\CCC$ (our $\ZZZ$), and computes $\deg \Gamma_{12} = 3762$ using the
genus 5 hyperelliptic locus $\FFF$  (see Section \ref{otherloci}).
He also discusses the loci $\BBB$ and $\DDD$.  (The quintic plane
curve with triple point from the proof of Prop. \ref{DA} appears
explicitly on p. 319.) He did not seem to be aware that the degree of
the morphisms $\CCC \rightarrow \Gamma_{12}$, $\DDD
\rightarrow \Gamma_{12}$ he describes are greater than 1.

\epoint{Coble}
Many of the constructions above have appeared, at least implicitly,
in \cite{co} (mainly Section 51, but also Sections 50 and 58).  

$\CCC$: On p. 210, the sextics with seven nodes from the proof of
Proposition \ref{CA} appear, and they are shown to correspond to the
120 pairs of (-1)-curves on a degree 1 Del Pezzo surface exchanged by
an involution, and also with odd theta-characteristics.  On p. 219,
this is connected to the data of plane quartics with another point in
the plane.

$\DDD$: On p.  108, the hyperelliptic curve of genus 3 appears, with
branch points corresponding to the 8 exceptionals on a degree 1 Del
Pezzo surface.  On p. 212, the quintic with a triple point from the
proof of Proposition \ref{DA} appears, although it is not identified
with the hyperelliptic genus 3 curve mentioned earlier.

$\EEE$: On p. 209, a construction from the proof of Proposition
\ref{EA} appears, the plane nonic with triple points at the 8 given
points.  The morphism $\AAA \rightarrow \EEE$ is explicitly described.
On p. 220, Coble remarks (section 58) that Schottky gives the coordinates of the 8
points and the equation of the nonic explicitly from modular functions
on the space $\cm^1_4$.

\section{Further questions}
\label{fq}

\epoint{Interpreting these results in terms of elliptic fibrations}
Given an elliptic fibration over some base $B$, the non-zero two-torsion
of the fibration is a triple cover of $B$ (and this can be tweaked to
give quadruple covers as 2-torsion information, see the proofs 
of Propositions \ref{CA} and \ref{DA}).  One natural question is:
what triple (or quadruple) covers can arise as two-torsion of
an elliptic fibration?

Also, if the family is not isotrivial, it can be reconstructed from the
$n$-torsion of the family if $n>2$.  A second natural question is:
to what extent is this true if $n=2$?

In the case of rational elliptic fibrations, both questions are 
answered completely.  For the first, striking explicit geometric conditions
are given.  And for the second, one can recover the fibration using the 
geometric conditions in the answer to the first.  

It would be interesting to extend this analysis to other 
elliptic fibrations.

(a) The next reasonable geometric example would be elliptic K3-surfaces, which
have 24 singular fibers.  This analog of $\AAA$ corresponds to a
codimension 3 locus in $\Sym^{24} \proj^1$, and is of course the same
as the analog of $\BBB$.  (Another codimension 3 locus in the
parameter space of 24 points on the projective line are genus 5 curves
mapped to $\proj^1$ by a canonical pencil, but there's no obvious reason
why this should be the same locus!)

(b)  Could any sense be made of this question arithmetically, e.g. over $\Q$?

} 

\end{document}